\documentclass{amsart}
\usepackage{amssymb}

\newcommand{\ZFCa}{{\operatorname{\mathsf {ZFC}}}}

\newcommand{\reals}{{\mathbb R}}

\newcommand{\rest}{{\mathord{\restriction}}}

\newcommand{\dom}{{\operatorname{\mathsf {dom}}}}
\newcommand{\ran}{{\operatorname{\mathsf    {range}}}}
\newcommand{\rank}{{\operatorname{\mathsf    {rank}}}}
\newcommand{\Seq}{{\operatorname{\mathsf    {Seq}}}}
\newcommand{\QED}{\hspace{0.1in} \square \vspace{0.1in}}
\newcommand{\forces}{\Vdash}
\newcommand{\V}{{\mathbf V}}
\newcommand{\HH}{{\mathbf H}}
\newcommand{\<}{\langle}
\renewcommand{\>}{\rangle}
\newcommand{\thinks}{\models}
\newcommand{\Proof}{{\sc Proof} \hspace{0.2in}}

\newcommand{\PPP}{{\mathcal P}}
\newcommand{\Ss}{\operatorname{\mathsf   {S}}}

\newcommand{\lft}[2]{\mathopen\ifcase#1{}\oo\or
                        \big#2\or\Big#2\else\oo\fi} 
\newcommand{\rgt}[2]{\mathclose\ifcase#1{}\oo\or
                        \big#2\or\Big#2\else\oo\fi}

\theoremstyle{plain}
\newtheorem{theorem}{Theorem}[section]
\theoremstyle{plain}
\newtheorem{lemma}[theorem]{Lemma}

\newtheorem{definition}[theorem]{Definition}

\begin{document}

\title{Splitting number}
\author{Tomek Bartoszy\'{n}ski}
\address{Department of Mathematics\\
Boise State University\\
Boise, Idaho 83725, USA}
\email{{\tt tomek@@math.idbsu.edu}}

\thanks{Research partially supported by 
NSF grant DMS 95-05375}
\maketitle

\begin{abstract}
We show that it is consistent with $\ZFCa$ that every uncountable set can be
continuously mapped onto a splitting family. 
  \end{abstract}
\section{Introduction}
A family $\mathcal A \subseteq [\omega]^\omega$ is called a splitting family if
for every infinite 
set $B \subseteq \omega$ there exists $A \in
\mathcal A$ such that 
$$|A \cap B|=|(\omega \setminus A) \cap B|=\boldsymbol\aleph_0.$$
We denote by $\mathfrak s$ the least size of a splitting family.
It is well-known that $ \boldsymbol\aleph_1 \leq \mathfrak s \leq  2^{
  \boldsymbol\aleph_0}$.
Let 
$$\Ss = \{X \subseteq 2^\omega: \text{ no Borel image of $X$ is a splitting
family}\}.$$
By ``Borel image'' we mean image by a Borel function.
It is easy to see that $\Ss$ is a $ \sigma $-ideal containing all
countable sets. The purpose of this paper is to show that one cannot
prove in $\ZFCa$ that $\Ss$ contains an uncountable set.

Recall (\cite{JuShSus88} or \cite{BJbook}) that a forcing notion
$(\PPP, \leq)$ is Suslin if 
\begin{enumerate}
  \item $\PPP$ is ccc,
  \item $\PPP$ is a ${\boldsymbol \Sigma}^1_1$ set of reals,
  \item relations $\leq, \perp$ are ${\boldsymbol \Sigma}^1_1$.
\end{enumerate}
Let ${\mathbf {MA}}(\text{Suslin})$ denote  Martin's Axiom for
Suslin partial orders.
It is well known that  ${\mathbf {MA}}(\text{Suslin})$ implies that
many cardinal invariants, most notably additivity of measure, are  equal to
$2^{\boldsymbol\aleph_0}$. 

Notation used in this paper is standard.
In particular, for $s,t \in  2^{<\omega}$, 
$[s]= \{x \in 2^\omega : x \rest \dom(s)=s\}$ and $s^\frown t$ denotes
the concatenation of $s$ and $t$. 
For $A,B \subseteq \omega$ define $A \subseteq^\star B$ if $|A
\setminus B| < \boldsymbol\aleph_0$.

To simplify notation,
throughout this paper we will identify elements of $[\omega]^\omega$ with
elements of $2^\omega$  via characteristic functions. 

\section{Consistency result}
The goal of this paper is to show ${\mathbf {MA}}(\text{Suslin})$ is
consistent with
$\Ss=[\reals]^{<\boldsymbol\aleph_1}$. 
This is a generalization of a result from \cite{JuShSus88}, where it
was proved that ${\mathbf {MA}}(\text{Suslin}) + \mathfrak s =
\boldsymbol\aleph_1$ is consistent.

\begin{theorem}
There exists a model of $\V' \thinks \ZFCa$ such that:
\begin{enumerate}
  \item $\V' \thinks \Ss=[\reals]^{<\boldsymbol\aleph_1}$,
  \item $\V' \thinks {\mathbf {MA}}(\text{Suslin}) +
    2^{\boldsymbol\aleph_0}=\boldsymbol\aleph_2$.
\end{enumerate}
\end{theorem}

The rest of this section is devoted to the proof of this theorem.

We start with a definition of a forcing notion, due to Hechler
(see \cite{he:cosub}),  that will be crucial
for our construction. 

Let $\Seq^\star \subseteq \omega^{<\omega}$ be the set of strictly
increasing finite sequences.  Let
$$ {\mathbf D}=\{ (s,f): s \in \Seq^\star, \ s \subseteq f, \ \text{ and
  $f$ is strictly increasing}\}.$$

For $(s,f),(t,g) \in \mathbf D$ define 
$$(s,f)  \geq (t,g) \iff
    s \supseteq t \ \& \  
\forall n \in \omega \ f(n) \geq g(n).$$

Define a rank function on ${\mathbf D}$:
\begin{definition}[\cite{BaumDor85Adj}]\label{drank}
Suppose that $D \subseteq {\mathbf D}$ is a
dense open set.
For $s \in \Seq^\star$ define the rank of $s$ as follows:
\begin{enumerate}
\item $\rank_{D}(s) = 0 $ if there exists a function $f$ such that
  $(s,f) \in D$.
\item If $\rank_{D}(s)\neq 0$, then
\begin{multline*}
\rank_{D}(s) = \min\lft2\{\alpha: \exists m \ 
    \exists \left\{s_k : k \in \omega\right\} \subseteq \Seq^\star
    \cap \omega^m \\
 \lft1(
 \rank_D(s_k) < \alpha \ \&\ s \subseteq
  s_k \ \&\ s_k(|s|)>k\rgt1)\rgt2\}.
\end{multline*}
\end{enumerate}
\end{definition}

\begin{lemma}[\cite{BaumDor85Adj}, \cite{BJbook} lemma 3.5.6]
  For every $s \in \Seq^\star$, $\rank_D(s)$ is defined.~$\QED$
\end{lemma}

Let $\<{\mathcal P}_\alpha, \dot{{\mathcal Q}}_\alpha:
\alpha<\omega_2)$ be a finite support iteration such that  
\begin{enumerate}
  \item $\forces_\alpha \dot{{\mathcal Q}}_\alpha \text{ is Suslin,}$
  \item if $\alpha$ is a limit ordinal then $\dot{{\mathcal
        Q}}_\alpha \simeq {\mathbf D}$.
\end{enumerate}
By careful bookkeeping we can ensure that $\V^{{\mathcal
    P}_{\omega_2}} \thinks {\mathbf {MA}}(\text{Suslin}) +
    2^{\boldsymbol\aleph_0}=\boldsymbol\aleph_2$.

We will show that $\V^{{\mathcal P}_{\omega_2}} \thinks
\Ss=[\reals]^{<\boldsymbol\aleph_1}$. The following construction is a
modification of a construction from \cite{BerJuCom}.

Suppose that $A \subseteq \omega$ is an infinite set. Let $A_-, A_+$
be two, canonically chosen, disjoint infinite sets such that $A_- \cup A_+ =A$.
Assume that $g \in \omega^\omega$ is an increasing function such that
$\ran(g) \cap A=\{x_n:n \in \omega\}$ is infinite ($x_n<x_{n+1}$ for
all $n$).

Define a real $z_{A,g} \in 2^\omega$ as follows:
$$z_{A,g}(n)=\left\{\begin{array}{ll}
1 & \text{if } x_n \in A_+\\
0 & \text{if } x_n \in A_-
\end{array}\right. \text{ for } n \in \omega .$$

Fix a bijection $\ell: 2^{<\omega} \longrightarrow \omega $ and for $x
\in 2^\omega$ define $\ell(x)=\{\ell(x \rest n): n \in \omega\}$.
Define 
$S_g: \dom(S_g) \longrightarrow 2^\omega $ as
$$S_g(x)=z_{\ell(x),g} \text{ for } x \in 2^\omega .$$

The following lists some easy properties of the function defined above:
\begin{lemma}
  \begin{enumerate}
  \item $\dom(S_g)$ is a  $G_\delta$ subset of $2^\omega $,
  \item $S_g$ is continuous on its domain,
    \item $S_g$ extends to a Borel function on $2^\omega $.
  \end{enumerate}
\end{lemma}
\Proof
(1) Note that $\dom(S_g) = \{x \in 2^\omega : |\ran(g) \cap \ell(x)|=
\boldsymbol\aleph_0\}$, which is a  $G_\delta $ set (possibly empty). 
(2) is easy to see and (3) is well-known.~$\QED$

\begin{definition}
  An uncountable set $X \subseteq 2^\omega $ is called a Luzin set if
  $|F \cap X| \leq \boldsymbol\aleph_0 $ for every meager set $F
  \subseteq 2^\omega $.
\end{definition}

We have the following easy lemma:
\begin{lemma}
  Every Luzin set is a splitting family.
\end{lemma}
\Proof 
Suppose that $X \subseteq 2^\omega $ is a non-meager  set. 
Since we identify elements of $[\omega]^\omega$ with
elements of $2^\omega$  via characteristic functions we can assume
that $X \subseteq [\omega]^\omega $. 
Let $A \in   [\omega]^\omega $.
Consider the set 
$$F=\{ z \in [\omega]^\omega : z \subseteq^\star ( \omega \setminus A)
\text{ or } A 
\subseteq^\star  z\}.$$
It is easy to see that $F$ is a meager set, and that  any element of $X
\setminus F$ splits $A$.~$\QED$

\begin{lemma}\label{main}
Suppose that $d$ is a ${\mathbf D}$-generic real over $\V$. 
If $Z \subseteq 2^\omega \cap \V$ is uncountable then $S_d(Z)$ is a
Luzin set in $\V[d]$.
\end{lemma}
\Proof
Observe first that by genericity $\V \cap 2^\omega \subseteq
\dom(S_d)$ and $S_d$ is one-to-one on $\V \cap 2^\omega$. In particular,
$S_d(Z)$ is an uncountable set.

Suppose that $F \in \V[d]$ is a closed nowhere dense subset of
$2^\omega$.  To
show that $S_d(Z)$ is a
Luzin set it is enough to show that $S_d(Z) \cap F$ is countable.
Let $f \in (2^{<\omega})^\omega \cap \V[d]$ be a function defined as
follows:

$$f(n)= \min\{s \in 2^{< \omega } : \forall t \in 2^{ \leq n}\ 
[t^\frown s] \cap F = \emptyset\}.$$
(The minimum is taken with respect to some canonical enumeration of $2^{<
  \omega }$.) It is well-known that such an $s$ exists.

Let $ \dot{f}$ be a $ {\mathbf D}$-name for $f$ and define for $ n \in
\omega $,
$$D_n =\{p \in {\mathbf D}: \exists s \in 2^{< \omega } \ p \forces_{{\mathbf D}}
\dot{f}(n)=s\}.$$ 

Let $N \prec \HH(\lambda)$ be a countable model containing $ \dot{f}$
and $Z$, where $ \lambda $ is a sufficiently large regular cardinal.
\begin{lemma}
If $x \in \V \cap 2^\omega$ but $x \not\in N \cap 2^\omega$ 
then $S_d(x) \not\in
F$.
\end{lemma}
\Proof
Suppose not and let $x \not\in N \cap 2^\omega$ be a counterexample.
Choose $(s,g)~\in~{\mathbf D}$ such that $$(s,g) \forces_{{\mathbf D}}
S_{\dot{d}}(x) \in \dot{F}.$$
Let $\widetilde{k}= |\ran(s) \cap \ell(x)|$. In other words, $S_d(x)
\rest \widetilde{k}$ is determined by $(s,g)$.
Let 
$$U = \{t \in \Seq^\star: s \subseteq t \ \& \ |\ran(t) \cap
\ell(x)|=\widetilde{k} \ \& 
\ \forall j \in \dom(t)\setminus \dom(s)\ t(j) \geq g(j)\}.$$

\begin{lemma}\label{4321a}
 $\min\{\rank_{D_{\tilde{k}}}(t): t \in U\}=0$. 
\end{lemma}
\Proof Suppose that the lemma is not true and let $t \in U$ be an element of
minimal rank.
By the definition there exists $m$ and a sequence
  $\left\{t_j : j \in \omega\right\} \subseteq \Seq^\star \cap \omega^m$
  such that for every $j$:
  \begin{enumerate}
  \item $t \subseteq t_j$,
  \item $\rank_{D_{\tilde{k}}}(t_j) <
  \rank_{D_{\tilde{k}}}(t)$,
\item $t_j(|t|)>j$.
  \end{enumerate}
Fix $i$ such that $|t| \leq i<m$ and let
$W_i = \{t_j(i) : j \in \omega\}$.
Note that every subsequence of $\{t_j : k \in \omega\}$ witnesses that
$\rank_{D_{\tilde{k}}}(t)>0$ as well.
Thus, by passing to a subsequence we can assume that there is 
a set $\ell(x_i)$ such
that $W_i \subseteq  \ell(x_i)$ or $W_i \cap \ell(x)$ is finite for
all $x$. In particular, if such a real $x_i$ exists it is a
member of $N$.

Since $x \not \in N$,  $W_i \cap \ell(x)$ is finite for all $|t| \leq i<m$.
Therefore, 
 there exists $j$ such that $\ran(t_j) \cap \ell(x) = \ran(t) \cap
\ell(x)$.
In particular, $t_j \in U$ and $\rank_{D_{\tilde{k}}}(t_j)
<\rank_{D_{\tilde{k}}}(t)$,
which is a contradiction.~$\QED$

Let $ t \in U$ be such that $\rank_{D_{\tilde{k}}}(t)=0$.
There exists $h \in \omega^\omega$ such that 
$(t,h) \in D_{\tilde{k}}$. Therefore, $(t , \max(h,g)) \geq (s,g)$ and
$(t , \max(h,g))$ decides the value of $\dot{f}(\widetilde{k})$.
Denote this value by 
$\widetilde{s}$. However, $(t , \max(h,g))$ does not put any
restrictions on values of 
$S_d(x)(j)$ for $j \geq \widetilde{k}$. Extend $t$ to $t'$ such that 
$$(t', \max(h,g)) \forces_{\mathbf D} \left(S_{\dot{d}}(x) \rest
\widetilde{k}\right)^\frown \widetilde{s} \subseteq S_{\dot{d}}(x).$$
It is clear that 
$$(t', \max(h,g)) \forces_{\mathbf D} S_{\dot{d}}(x) \not\in \dot{F}.$$
This contradiction ends the proof of lemma \ref{main}.~$\QED$

Let $G \subseteq {\mathcal P}_{\omega_2}$ be a generic filter over $\V$.
Suppose that $Z \subseteq 2^\omega \cap \V[G]$ is a set of cardinality
$ \boldsymbol\aleph_1 $. First we find a limit ordinal 
$\alpha$ such that $Z \subseteq
\V[G \cap {\mathcal P}_\alpha]$. We will work in the model
$\V_1=\V[G \cap {\mathcal P}_{\alpha+1}] = \V[G \cap {\mathcal
  P}_\alpha][d]$, where $d$ is a ${\mathbf D}$-generic real over $\V[G
\cap {\mathcal P}_\alpha]$.

To finish the proof it is enough to show that $S_d(Z)$ is a splitting
family in $\V[G]$. Note however that $S_d(Z)$ is not a Luzin set in
$\V[G]$. In fact, $S_d(Z)$ is meager in
$\V[G]$.

\begin{lemma}
 $\{S_d(x): x \in Z\}$ is a splitting
family in $\V[G]$.
\end{lemma}
\Proof We will work in $\V_1$.
By \ref{main}, we know that $S_d(Z)$ is a Luzin  set in
 in $\V_1$. 
Note that $\V[G]$ is a generic extension of $\V_1$ via finite support
iteration of Suslin forcings $ {\mathcal P}_{\alpha+1, \omega_2}$.

Let $\dot{A}$ be a
$\PPP_{\alpha+1,\omega_2}$-name for a set $A \in [\omega]^\omega$. 
We will need the following lemma:
\begin{lemma}
  For every $p \in {\mathcal P}_{\alpha+1, \omega_2}$, the set
$$Z_p=\left\{z \in Z: p \forces_{\alpha+1, \omega_2} \text{``}S_d(z)
  \subseteq^\star (\omega \setminus 
\dot{A})   \text{ or }  \dot{A}
\subseteq^\star S_d(z)\text{''}\right\}$$
is countable.
\end{lemma}
Before we prove the lemma notice that the theorem follows from it
immediately -- given $p \in \PPP_{\alpha+1,\omega_2}$, $\dot{A}$ and $z \in Z
\setminus Z_p$ we can find $q \geq p$ such that 
$q \forces \text{``$S_d(z)$ splits $ \dot{A}$.''}$

{\sc Proof of the lemma} 
We will use the absoluteness properties of Suslin forcing (see
\cite{BJbook} or \cite{JuShSus88}).

Fix a condition $p \in {\mathcal P}_{\alpha+1, \omega_2}$. 
Let
$M$ be a countable 
elementary submodel of $\HH(\lambda)$ containing $\dot{A}$, $p$ and
$\PPP_{\alpha+1,\omega_2}$. 
Define a finite support iteration $\<\PPP_\alpha(M), \dot{{\mathcal
    Q}}_\alpha(M): \alpha<\omega_2\>$ as follows:

$$\forces_\alpha \dot{{\mathcal Q}}_\alpha(M) = \left\{\begin{array}{ll}
\dot{{\mathcal Q}}_\alpha & \text{if $\alpha \in M$}\\
\emptyset & \text{if $\alpha \not\in M$}
\end{array}\right. \text{ for $\alpha<\omega_2$}. $$
Let ${\mathcal P}=\lim {\mathcal P}_\alpha(M)$. $ {\mathcal P} $ is
the part of the iteration that
contains all information regarding
$ \dot{A}$.
${\mathcal P}$ is isomorphic to a countable iteration of Suslin
forcings. In particular, ${\mathcal P}$ has a definition that
can be coded as a real number (essentially by encoding $M$ as a real number). 

{}From Suslinness it follows  that ${\mathcal P} \lessdot {\mathcal
  P}_{\alpha+1,\omega_2}$ and 
that $\dot{A}$ is a ${\mathcal P}$-name (see \cite{BJbook} lemma 9.7.4 or
\cite{JuShSus88}). Moreover,
it is enough to show that
$$\left\{z \in Z: p \forces_{{\mathcal P}} \text{``}S_d(z)
  \subseteq^\star (\omega \setminus 
\dot{A})   \text{ or }  \dot{A}
\subseteq^\star S_d(z)\text{''}\right\}$$
is countable.

Let $N \prec \HH(\lambda)$ be a countable model containing $M, \dot{A}$
and ${\mathcal P}$.
Since $S_d(Z)$ is a Luzin set in $\V_1$, the set 
$$Z_0=\{z \in Z: S_d(z) \text{ is not a Cohen real over } N\}$$
is countable. 
We will show that $Z_p \subseteq Z_0$. In particular, 
for $z \in Z \setminus Z_0$, 
$$p \not \forces_{\PPP}
\text{``}S_d(z) \subseteq^\star (\omega \setminus \dot{A}) \ \text{ or
  } \ \dot{A}  
\subseteq^\star S_d(z)\text{''},$$
 which will finish the proof.  
Fix 
$z \in Z \setminus Z_0$ and let $Y=S_d(z)$ be a Cohen real over $N$.
Without loss of generality we can assume  that $p \forces_{{\mathcal
    P}} Y \subseteq^\star (\omega \setminus \dot{A})$. 

Clearly, $N[Y][G \cap N[Y]] \thinks Y \subseteq^\star (\omega
\setminus \dot{A}[G \cap N[Y]])$
and therefore  
$$N[Y] \thinks \text{``$p \forces_{{\mathcal P}} Y \subseteq^\star
  (\omega \setminus \dot{A})$,''}$$
since the last statement is absolute.
Represent the Cohen algebra as $ \mathbf C = [\omega]^{< \omega}$ and let
$ \dot{Y}$ be the canonical name for a Cohen real.
There is a condition $s \in \mathbf C$ such that 
$$N \thinks s \forces_{\mathbf C} \text{``}p \forces_{ {\mathcal P}} \dot{Y}
\subseteq^\star (\omega \setminus \dot{A}).\text{''}$$

Let $Y'= s \cup \lft2(\omega \setminus \lft1(Y \setminus
\max(s)\rgt1)\rgt2)$.
$Y'$ is also a a Cohen real over $N$ and since $s \subseteq Y'$ we get
that
$N[Y'] \thinks \text{``$p \forces_{{\mathcal P}} Y' \subseteq^\star
  (\omega \setminus \dot{A})$.''}$
It follows that 
$$N[Y'][G \cap N[Y']] \thinks Y' \subseteq^\star (\omega \setminus \dot{A}[G
\cap N[Y']]).$$ 
Note that $ \dot{A}[G] = \dot{A}[G
\cap N[Y']]=\dot{A}[G \cap N[Y]]$.
Thus $\V[G] \thinks Y \cup Y' \subseteq^\star (\omega \setminus \dot{A}[G])$ which means that $
\dot{A}[G]$ is finite. Contradiction.

The same argument shows that the assumption that $p \forces_{{\mathcal
    P}}   \dot{A} \subseteq^\star Y $ leads to a contradiction.~$\QED$

{\bf Acknowledgement}: I would like to thank Andreas Blass for his
helpful communication concerning the preparation of this paper. 

\ifx\undefined\bysame
\newcommand{\bysame}{\leavevmode\hbox to3em{\hrulefill}\,}
\fi


\begin{thebibliography}{1}

\bibitem{BJbook}
Tomek Bartoszy\'{n}ski and Haim Judah, {\em Set {T}heory: on the structure of
  the real line}, A.K. Peters, 1995.

\bibitem{BaumDor85Adj}
James~E. Baumgartner and Peter Dordal, {\em Adjoining dominating functions},
  The Journal of Symbolic Logic {\bf 50} (1985), no.~1, 94--101.

\bibitem{BerJuCom}
Jorg Brendle, Haim Judah, and Saharon Shelah, {\em Combinatorial properties of
  {H}echler forcing}, Annals of Pure and Applied Logic {\bf 59} (1992),
  185--199.

\bibitem{he:cosub}
S.~H. Hechler, {\em On the existence of certain cofinal subsets of
  $\omega^\omega$}, Axiomatic Set Theory (T.~J. Jech, ed.), Proc. Symp. Pure
  Math., vol.~13, Amer. Math. Soc., Providence, R.I., 1974, Part 2,
  pp.~155--173.

\bibitem{JuShSus88}
Haim Judah and Saharon Shelah, {\em Suslin forcing}, The Journal of Symbolic
  Logic {\bf 53} (1988), 1188--1207.

\end{thebibliography}

\end{document}